\documentclass[a4paper,10pt,reqno, english]{amsart} 
\usepackage{amssymb} 
\usepackage{amsmath,amsthm}
\usepackage{amsfonts,amssymb,enumerate}
\usepackage{url,paralist}
\usepackage{mathtools} 
\usepackage{anysize} 
\usepackage{hyperref} 
\usepackage{enumerate,amssymb}
\usepackage{mathtools}
 
\theoremstyle{plain}
\newtheorem{theorem}{Theorem}[section]

\theoremstyle{definition}

\newcommand{\BIGOP}[1]{\mathop{\mathchoice%
{\raise-0.22em\hbox{\huge $#1$}}%
{\raise-0.05em\hbox{\Large $#1$}}{\hbox{\large $#1$}}{#1}}}

\newcommand{\R}{\mathbb{R}} 
\newcommand{\N}{\mathbb{N}}
\newcommand{\Z}{\mathbb{Z}}
\newcommand{\F}{\mathbb{F}} 
\newcommand\Sym{\mathfrak S} 
\newcommand{\ind}{\operatorname{Index}}

\begin{document}

\title{Polynomial partitioning for several sets of varieties}

\author[Blagojevi\'c]{Pavle V. M. Blagojevi\'{c}} 
\thanks{The research by Pavle V. M. Blagojevi\'{c} leading to these results has
        received funding from DFG via the Collaborative Research Center TRR~109 ``Discretization in Geometry and Dynamics.''.
        Also supported by the grant ON 174008 of the Serbian Ministry of Education and Science.}
\address{Inst. Math., FU Berlin, Arnimallee 2, 14195 Berlin, Germany\hfill\break%
\mbox{\hspace{4mm}}Mat. Institut SANU, Knez Mihailova 36, 11001 Beograd, Serbia}
\email{blagojevic@math.fu-berlin.de} 
\author[Dimitrijevi\'c Blagojevi\'c]{Aleksandra S. Dimitrijevi\'c Blagojevi\'c}
\thanks{The research by Aleksandra Dimitrijevi\'c Blagojevi\'c leading to these results has
        received funding from the  grant ON 174008 of the Serbian Ministry of Education and Science.}
\address{Mat. Institut SANU, Knez Mihailova 36, 11001 Beograd, Serbia}
\email{aleksandra1973@gmail.com}
\author[Ziegler]{G\"unter M. Ziegler} 
\thanks{The research by G\"unter M. Ziegler received funding from DFG via the Research Training Group “Methods for Discrete Structures” and the Collaborative Research Center TRR~109 ``Discretization in Geometry and Dynamics.''}  
\address{Inst. Math., FU Berlin, Arnimallee 2, 14195 Berlin, Germany} 
\email{ziegler@math.fu-berlin.de}


\begin{abstract}
We give a new, systematic proof for a recent result of Larry Guth and thus also extend the result to a setting with several families of varieties:
For any integer $D\geq 1$ and any collection of sets $\Gamma_1,\ldots,\Gamma_j$ of low-degree $k$-dimensional varieties in~$\R^n$ there exists a non-zero polynomial $p\in\R[X_1,\ldots,X_n]$ of degree at most $D$ so that each connected component of $\R^n{\setminus}Z(p)$ intersects $O(jD^{k-n}|\Gamma_i|)$ varieties of $\Gamma_i$, simultaneously for every $1\leq i\leq j$.
For $j=1$ we recover the original result by Guth.
Our proof, via an index calculation in equivariant cohomology, 
shows how the degrees of the polynomials used for partitioning are dictated by the topology, 
namely by the Euler class being given in terms of a top Dickson polynomial.
\end{abstract}


\date{}

\maketitle


\section{Introduction}
\label{sec : Introduction}

The celebrated work \cite{guth2015-2} by Larry Guth and Nets Hawk Katz
on the Erd\H{o}s distinct distances problem 
in the plane brought to light the following beautiful partitioning result:
\begin{theorem}[Guth and Katz 2015 {\cite[Thm.\,4.1]{guth2015-2}}]
\label{th : guth-katz}
	Let $X$ be a finite set of points in $\R^n$, and let $D\geq 1$ be an integer.
	Then there exists a non-zero polynomial $p\in\R[X_1,\ldots,X_n]$ of degree at most $D$ such that each connected component of the complement $\R^n{\setminus}Z(p)$ contains at most $C_nD^{-n}|X|$ points of $X$, where $C_n$ is a constant that may depend on $n$.
\end{theorem}

\noindent
Here $Z(p)$ denotes the set of zeroes in $\R^n$ of the polynomial $p$, that is 
\[
Z(p)=\{ (x_1,\ldots,x_n)\in\R^n : p(x_1,\ldots,x_n)=0\}.
\]

In his recent paper \cite{guth2015}, Guth used  equivariant topology to prove the following extended polynomial partitioning result. 

\begin{theorem}[Guth,\,2015 {\cite[Thm.\,0.3]{guth2015}}]
\label{th : guth}
	Let $\Gamma$ be a finite set of $k$-dimensional varieties in $\R^n$, each of them defined by at most $m$ polynomial equations of degree at most $d$.
	Then for any $D\geq 1$ there exists a non-zero polynomial $p\in\R[X_1,\ldots,X_n]$ of degree at most $D$ such that each connected component of the complement $\R^n{\setminus}Z(p)$ intersects at most $C(d,m,n)D^{k-n}|\Gamma|$ varieties in $\Gamma$, where $C(d,m,n)$ is a constant that may depend on the parameters $d$, $m$, and $n$.
\end{theorem}
 
In this paper, based on the set-up from the proof of the previous theorem and the use of the Fadell--Husseini index theory \cite{fadell1988} for the proof of a necessary Borsuk--Ulam type theorem, we make the next extension step by proving the following  ``colored'' generalization of Theorem~\ref{th : guth}.

\begin{theorem}
\label{th : main}
	Let $j\geq 1$ be an integer. 
	For $1\leq i\leq j$, let $\Gamma_i$ be a finite set of $k_i$-dimensional varieties in $\R^n$, each of them defined by at most $m_i$ polynomial equations of degree at most $d_i$.
	Then for any $D\geq 1$ there exists a non-zero polynomial $p\in\R[X_1,\ldots,X_n]$ of degree at most $D$ such that each connected component of the complement $\R^n{\setminus}Z(p)$ for every $1\leq i\leq j$ intersects at most $C(d_i,m_i,n)jD^{k_i-n}|\Gamma_i|$ varieties in~$\Gamma_i$, where $C(d_i,m_i,n)$ is a constant that may depend on parameters $d_i$, $m_i$, and $n$.
\end{theorem}

In a concrete example, this says the following:
There are constants $C_1=C(1,2,3)$ and $C_2=C(1,3,3)$
such that if we have (large) collections $\Gamma_1$ of red lines and $\Gamma_2$ of blue points in~$\R^3$,
then for every $D\ge1$ there is a nonzero polynomial $p(x,y,z)\in\R[x,y,z]$ of degree at most~$D$
such that each connected component of $\R^3\setminus Z(p)$ meets 
	at most $2C_1\dfrac{|\Gamma_1|}{D^2}$ red lines, and
	at most $2C_2\dfrac{|\Gamma_2|}{D^3}$ blue points.  
(This is the special case when we have $j=2$ families of varieties in $\R^3$, so $n=3$, where
the first family consists of lines, so $k_1=1$ and e.g.\ $m_1=2$, $d_1=1$,
and the second one of points, so $k_2=0$, $m_2=3$, $d_2=1$.)

\noindent
\emph{Note added in revision.}
In his recent paper \cite{BenYang_Joints}, Ben Yang used Theorem \ref{th : main} to obtain new upper bunds on the number of generalized joints for (possibly higher-dimensional) varieties.

\section{Proof of Theorem~\ref{th : main}}
\label{sec : proof}

The proof of Theorem~\ref{th : main} will have several separate components that at the end of the proof merge into the final argument.
We also rely on two particular results by J\'ozsef Solymosi and Terence Tao \cite[Thm.\,A.2]{solymosi2012} and by Guth \cite[Lemma\,3.1]{guth2015}.

\subsection{}
Let $P^{\delta}_n$ be the vector space of polynomials in $n$ variables of degree at most $\delta$
 with real coefficients, 
\[
P^{\delta}_n=\{ p\in \R[X_1,\ldots,X_n] : \deg p \leq\delta \}.
\]
The dimension of this vector space is $\dim P^{\delta}_n= {\delta +n \choose n}>\tfrac{{\delta}^n}{n!}$.
For every integer $\ell\geq 1$ choose the smallest integer $\delta_{\ell}$ with the property that
\begin{equation*}
	j2^{\ell-1}\leq \tfrac{{\delta_{\ell}^n}}{n!}<  j2^n2^{\ell-1},
\end{equation*}
or equivalently 
\begin{equation}
	\label{ineq : 01}
	(n!)^{\frac1n}j^{\frac1n}2^{\frac{\ell-1}{n}}\leq \delta_{\ell}< 2(n!)^{\frac1n}j^{\frac1n}2^{\frac{\ell-1}{n}}.
\end{equation}
In particular, $\dim P^{\delta_{\ell}}_n>j2^{\ell-1}$. 
Next, let $s$ be the smallest integer such that
\begin{equation}
	\label{ineq : 02}
	\sum_{\ell=1}^s\delta_{\ell}\leq D< \sum_{\ell=1}^{s+1}\delta_{\ell}.
\end{equation}
The inequalities \eqref{ineq : 01} and \eqref{ineq : 02} imply that
\begin{equation}
\label{ineq : 03}
D< \sum_{\ell=1}^{s+1}\delta_{\ell} < 2 (n!)^{\frac1n}j^{\frac1n} \sum_{\ell=1}^{s+1} 2^{\frac{\ell-1}{n}}=
2 (n!)^{\frac1n}j^{\frac1n}\, \frac{2^{\frac{s+1}{n}}-1}{2^{\frac{1}{n}}-1}
<\frac{2 (n!)^{\frac1n}}{2^{\frac{1}{n}}-1} \, j^{\frac1n}2^{\frac{s+1}{n}}
= \frac{2^{1+\frac1n} (n!)^{\frac1n}}{2^{\frac{1}{n}}-1} \, j^{\frac1n}2^{\frac{s}{n}}.
\end{equation}
Consequently, the inequality \eqref{ineq : 03} gives
\begin{equation}
\label{ineq : 04}
D^n< \frac{2^{n+1} n!}{(2^{\frac{1}{n}}-1)^n} \, j2^{s}
\qquad\Longrightarrow\qquad
\frac{1}{2^{s}}< \frac{2^{n+1} n!}{(2^{\frac{1}{n}}-1)^n}\, \frac{j}{D^n} = C_n  \frac{j}{D^n},
\end{equation}
where $C_n$ depends only on $n$.

\subsection{}
For every $1\leq\ell\leq s$ we have that $\dim P^{\delta_{\ell}}_n\geq j2^{\ell-1}+1$.
Let $V_{\ell}$ denote an arbitrary vector subspace of $P^{\delta_{\ell}}_n$ of dimension $j2^{\ell-1}+1$.
The unit sphere $S(V_{\ell})$ in the vector space $V_{\ell}$ is equipped with the free $\Z/2=\langle\omega_{\ell}\rangle$ action given by $\omega_{\ell} \cdot p=-p$ for $p \in S(V_{\ell})$.

\medskip
We will use the product space 
\begin{equation*}
	Y:=\prod_{\ell=1}^s S(V_{\ell})\cong \prod_{\ell=1}^s S^{j2^{\ell-1}}.
\end{equation*}
The elementary abelian group $(\Z/2)^s=\langle \omega_1,\ldots,\omega_s\rangle$ acts on $Y$ componentwise, that is, for $1\leq\ell\leq s$ and $(p_1,\ldots,p_{\ell},\ldots,p_s)\in Y$ the generator $\omega_{\ell}$ acts as follows:
\begin{equation}
	\label{eq : Y}
	\omega_{\ell}\cdot (p_1,\ldots,p_{\ell},\ldots,p_s)=(p_1,\ldots,-p_{\ell},\ldots,p_s).
\end{equation}

\subsection{}
Consider the vector space $\R^{(\Z/2)^s}$ and the vector subspace of codimension $1$ given by
\begin{equation*}
U_s=\Big\{ (y_{\alpha})_{\alpha\in (\Z/2)^s} \in \R^{(\Z/2)^s} : \sum_{\alpha\in (\Z/2)^s} y_{\alpha} =0\Big\}.
\end{equation*}
We introduce the following action of $(\Z/2)^s$ on $\R^{(\Z/2)^s}$:
The element $(\beta_1, \ldots, \beta_s) \in (\Z/2)^s$  acts on the vector $(y_{\alpha})_{\alpha \in (\Z/2)^s}\in \R^{(\Z/2)^s}$ by acting on its index set
\begin{equation}
	\label{eq : U_s}
	(\beta_1, \dots, \beta_s)  \cdot (\alpha_1, \dots, \alpha_s) = (\beta_1 + \alpha_1, \ldots, \beta_s + \alpha_s),
\end{equation}
where the addition is assumed to be in $\Z/2$. 
With respect to the introduced action the vector subspace $U_s$ is a $(\Z/2)^s$-subrepresentation of $\R^{(\Z/2)^s}$ of dimension $2^s-1$.
%

\subsection{}

Any non-constant polynomial $p\in \R[X_1,\ldots,X_n]$ determines two disjoint open regions
in~$\R^n$, possibly one of them empty, which we denote by 
\[
D_p^0=\{(x_1,\ldots,x_n)\in \R^n : p(x_1,\ldots,x_n)>0\}\quad\text{and}\quad
D_p^1=\{(x_1,\ldots,x_n)\in\R^n : p(x_1,\ldots,x_n)<0\}.
\]
Thus $\R^n{\setminus} Z(p)=D_p^0\cup D_p^1$ and $D_p^0\cap D_p^1 =\emptyset$.

\medskip
Let $(p_1,\ldots,p_s)\in Y$ be an ordered tuple of polynomials in $Y$, and let $\alpha=(\alpha_1,\ldots,\alpha_s)\in (\Z/2)^s=\{0,1\}^s$.
The \emph{sign pattern domain} determined by the tuple $(p_1,\ldots,p_s)$ and by the element $\alpha\in(\Z/2)^s$ is the intersection of open regions 
\[
\mathcal{O}_{\alpha}^{(p_1,\ldots,p_s)}=D_{p_1}^{\alpha_{1}}\cap\cdots\cap D_{p_s}^{\alpha_{s}}.
\]
An sign pattern domain can be empty. 
Moreover
\begin{equation}
	\label{eq : complement}
	\R^n{\setminus}Z(p_1\cdots p_s)=\bigcup_{\alpha\in (\Z/2)^s} \mathcal{O}_{\alpha}^{(p_1,\ldots,p_s)},
\end{equation}
where the union is disjoint union.
Observe that $\deg (p_1\cdots p_s)\leq\sum_{\ell=0}^s\delta_\ell \leq D$. 
Furthermore, the sign pattern domains $\mathcal{O}_{\alpha}^{(p_1,\ldots,p_s)}$ are unions of connected components of the complement $\R^n{\setminus}Z(p_1\cdots p_s)$.
\subsection{}
\label{subsec : 2.5}
 
For every $\alpha\in (\Z/2)^s$ and every variety $\gamma\subset\R^n$ we define the function $\phi_{\alpha,\gamma}\colon Y\longrightarrow\R$ by
\[
\phi_{\alpha,\gamma} (p_1,\ldots,p_s)=
\begin{cases}
	1, & \text{if } \mathcal{O}_{\alpha}^{(p_1,\ldots,p_s)}\cap \gamma\neq\emptyset\\
	0, & \text{if } \mathcal{O}_{\alpha}^{(p_1,\ldots,p_s)}\cap \gamma = \emptyset,
\end{cases}
\]
where $(p_1,\ldots,p_s)\in Y$.  
The functions $\phi_{\alpha,\gamma}$ are not continuous, but as Guth showed in \cite[Lemma\,3.1]{guth2015}, they can be approximated by sequences of continuous functions:
\begin{quote}
	{\small
	{\bf Lemma 3.1.} 
	{\em For $\varepsilon>0$, $\gamma\subset\R^n$ and $\alpha\in  (\Z/2)^s$, we define functions $\phi_{\alpha,\gamma}^{\varepsilon}\colon Y\longrightarrow\R$ with the following properties.
	\begin{compactenum}[\rm (1)]
		\item The functions $\phi_{\alpha,\gamma}^{\varepsilon}\colon Y\longrightarrow\R$ are continuous.
		\item $0\leq \phi_{\alpha,\gamma}^{\varepsilon} \leq 1$.
		\item If $\mathcal{O}_{\alpha}^{(p_1,\ldots,p_s)}\cap \gamma = \emptyset$, then $\phi_{\alpha,\gamma}^{\varepsilon}=0$.
		\item If $\varepsilon_i\rightarrow 0$ and $(p_1^i,\ldots,p_s^i)\rightarrow (p_1,\ldots,p_s)$ in $Y$, and $\mathcal{O}_{\alpha}^{(p_1,\ldots,p_s)}\cap \gamma\neq\emptyset$, then 
			\[
			\lim_{i\rightarrow\infty} \phi_{\alpha,\gamma}^{\varepsilon_i}(p_1^i,\ldots,p_s^i)=1.
			\]
		In other words, $\phi_{\alpha,\gamma} (p_1,\ldots,p_s)\leq\liminf_{i\rightarrow\infty}	\phi_{\alpha,\gamma}^{\varepsilon_i}(p_1^i,\ldots,p_s^i)$.
	\end{compactenum}
	}}
\end{quote} 
In order to simplify the presentation we postpone the typical compactness argument applied to $Y$ to the last step in the proof 
(\ref{subsec:compactness_argument}) and for now continue to work with the functions $\phi_{\alpha,\gamma}$ as if they were continuous.

\subsection{}

Let $1\leq i\leq j$ be fixed.
By assumption $\Gamma_i$ is a finite set of $k_i$-dimensional varieties in $\R^n$, each defined by at most $m_i$ polynomial equations of degree at most $d_i$.
Consider the following map from the space $Y$ to the representation $U_{s}^{\oplus j}$ associated to the collection of finite sets of varieties $\Gamma_1,\ldots,\Gamma_j$ given in the theorem:
\begin{align}
\label{eq : def of Phi}
\Phi \colon Y &\longrightarrow U_{s}^{\oplus j},\nonumber \\
(p_1,\ldots,p_s) &\longmapsto \Big( \Big( 
\sum_{\gamma\in\Gamma_i}\phi_{\alpha,\gamma}(p_1,\ldots,p_s)-\dfrac{1}{2^s}
\sum_{\beta \in (\Z/2)^s} \sum_{\gamma\in\Gamma_i}\phi_{\beta,\gamma}(p_1,\ldots,p_s)
\Big)_{\alpha\in (\Z/2)^s} \Big)_{i \in \{1,\dots, j\}}.
\end{align}
The sum $\sum_{\gamma\in\Gamma_i}\phi_{\alpha,\gamma}(p_1,\ldots,p_s)$ 
counts the number of varieties in $\Gamma_i$ that intersect sign pattern domain $\mathcal{O}_{\alpha}^{(p_1,\ldots,p_s)}$.
The map $\Phi$ is continuous and $(\Z/2)^s$-equivariant with respect to the actions given by \eqref{eq : Y} and~\eqref{eq : U_s}, 
assuming the diagonal action on the direct sum $U_{s}^{\oplus j}$.

\medskip
Let us assume that $\Phi^{-1}(0)\neq\emptyset$,
and pick $(p_1,\ldots,p_s)\in \Phi^{-1}(0)$.
Then each of $2^s$ sign pattern domains $\mathcal{O}_{\alpha}^{(p_1,\ldots,p_s)}$, $\alpha\in(\Z/2)^s$, determined by the tuple $(p_1,\ldots,p_s)$ in $Y$ intersects the same number of varieties in the set $\Gamma_i$, for every $1\leq i\leq j$. 
We use the following result of Solymosi and Tao \cite[Thm.\,A.2]{solymosi2012}, stated as in \cite[Thm.\,0.2]{guth2015}:
\begin{quote}
	{\small
	{\bf Theorem A.2.} 
	{\em  
	Suppose that $\gamma$ is a $k$-dimensional variety in $\R^n$ defined by $m$ polynomial equations each of degree at most $d$.
	If $p$ is a polynomial of degree at most $D$, then $\gamma$ intersects at most $C'(d,m,n)D^k$ different connected components of $\R^n{\setminus}Z(p)$, where $C'(d,m,n)$ is a constant that may depend on the parameters $d_i$, $m_i$, and $n$.
		}}
\end{quote}
Hence, each variety $\gamma\in\Gamma_i$ intersects at most $C'(d_i,m_i,n)D^{k_i}$ 
connected components of the complement $\R^n{\setminus}Z(p_1\cdots p_s)$.
Since each sign pattern domain $\mathcal{O}_{\alpha}^{(p_1,\ldots,p_s)}$ is a disjoint union of 
connected components of  $\R^n{\setminus}Z(p_1\cdots p_s)$, each variety $\gamma\in\Gamma_i$ 
intersects at most $C'(d_i,m_i,n)D^{k_i}$ sign pattern domains.
As we are looking at a point $(p_1,\ldots,p_s)\in \Phi^{-1}(0)$,
where we get the same number of varieties $\gamma\in\Gamma_i$ intersecting each sign pattern domain $\mathcal{O}_{\alpha}^{(p_1,\ldots,p_s)}$, this number is at most
$\tfrac1{2^s}|\Gamma_i|C'(d_i,m_i,n)D^{k_i}$.
The inequality \eqref{ineq : 04} implies that 
\[
\frac1{2^s}|\Gamma_i|C'(d_i,m_i,n)D^{k_i} < C_nC'(d_i,m_i,n) jD^{k_i-n}|\Gamma_i|.
\]
Each connected component of $\R^n{\setminus}Z(p_1\cdots p_s)$ is contained in a unique sign pattern domain, and therefore the number of varieties $\gamma\in\Gamma_i$ intersecting a connected component of $\R^n{\setminus}Z(p_1\cdots p_s)$ cannot exceed 
\[
C(d_i,m_i,n)\cdot jD^{k_i-n}|\Gamma_i|,
\]
where $C(d_i,m_i,n)$ is a constant that may depend on the parameters $d_i$, $m_i$, and $n$. 
This concludes the proof of Theorem~\ref{th : main}, except that it remains
to be verified that the map $\Phi$ indeed has a zero.

\subsection{}
We still need to prove that the $(\Z/2)^s$-equivariant map $\Phi \colon Y \longrightarrow U_{s}^{\oplus j}$, 
defined in \eqref{eq : def of Phi}, has a zero. Indeed, 
in the spirit of the usual resolution of ``configuration space/test map schemes'' for discrete geometry problems
\cite{Matousek:BU},
we will show that there is no continuous $(\Z/2)^s$-equivariant map $Y \longrightarrow U_{s}^{\oplus j}$ at all that avoids zero.

\medskip
Let us the assume to the contrary that there is such a map, then 
this induces a continuous $(\Z/2)^s$-equivariant map $Y \longrightarrow S(U_{s}^{\oplus j})$, 
where $S(U_{s}^{\oplus j})$ denotes the unit sphere in the vector spaces $U_{s}^{\oplus j}$.
Using the Fadell--Husseini ideal-valued index theory \cite{fadell1988} for the group $(\Z/2)^s$ and $\F_2$ coefficients, we will prove that such an equivariant map cannot exist, obtaining the required contradiction.

The cohomology of the group $(\Z/2)^s$ with $\F_2$ coefficients is given by $H^*((\Z/2)^s;\F_2)\cong\F_2[u_1,\ldots,u_s]$, where $\deg(u_i)=1$ and the variable $u_i$ corresponds to the generator $\omega_i$, for $1\leq i\leq s$.
According to \cite[Ex.\,3.3]{fadell1988},
\[
\ind_{(\Z/2)^s} (Y;\F_2)=\big\langle u_1^{j+1},u_2^{2j+1},u_3^{4j+1}\ldots, u_s^{j2^{s-1}+1}\big\rangle.
\]
Furthermore, from \cite[Prop.\,3.13]{fadell1988} (see also \cite[Prop.\,3.12 and Prop.\,3.13]{blagojevic2011}) 
we have that 
\[
\ind_{(\Z/2)^s}(S(U_{s}^{\oplus j});\F_2)=\Big\langle  \Big(\prod_{(\alpha_1,\ldots,\alpha_s)\in (\Z/2)^s{\setminus}\{0\}}  (\alpha_1u_1+\cdots+\alpha_su_s)\Big)^j \Big\rangle.
\]
Since a continuous $(\Z/2)^s$-equivariant map $Y \longrightarrow S(U_{s}^{\oplus j})$ exists, the basic property of the Fadell--Husseini index \cite[Sec.\,2]{fadell1988} implies that 
\[
\ind_{(\Z/2)^s}(S(U_{s}^{\oplus j});\F_2)\subseteq \ind_{(\Z/2)^s} (Y;\F_2), 
\]
and consequently that
\begin{equation}
\label{eq:upper_bound}
\Big(\prod_{(\alpha_1,\ldots,\alpha_s)\in (\Z/2)^s{\setminus}\{0\}}  (\alpha_1u_1+\cdots+\alpha_su_s)\Big)^j\in\big\langle u_1^{j+1},u_2^{2j+1},u_3^{4j+1}\ldots, u_s^{j2^{s-1}+1}\big\rangle.
\end{equation}
The polynomial
\[
q:=\prod_{(\alpha_1,\ldots,\alpha_s)\in (\Z/2)^s{\setminus}\{0\}}  (\alpha_1u_1+\cdots+\alpha_su_s)\,\in\,\F_2[u_1,\dots,u_s].
\]
is the Dickson polynomial of maximal degree \cite[Sec.\,III.2]{adem_milgram2004}.
It can be presented in the form
\[ q=\sum_{\pi\in\Sym_s}u_{\pi(1)}^{2^{s-1}}u_{\pi(2)}^{2^{s-2}}\cdots u_{\pi(s)}^{2^0}.\] 
Now the $j$-th power of the Dickson polynomial $q^j$ can be decomposed as follows
\begin{eqnarray*}
q^j &=& \Big(\prod_{(\alpha_1,\ldots,\alpha_s)\in (\Z/2)^s{\setminus}\{0\}}  (\alpha_1u_1+\cdots+\alpha_su_s)\;\; \Big)^j
=
\Big(\sum_{\pi\in\Sym_s}u_{\pi(1)}^{2^{s-1}}u_{\pi(2)}^{2^{s-2}}\cdots u_{\pi(s)}^{2^0}\Big)^{j}\\
       &=&\big(u_{s}^{j2^{s-1}}u_{s-1}^{j2^{s-2}}\cdots u_{2}^{2j} u_{1}^{j}\big) +\mathrm{Rest},
\end{eqnarray*}
where ``$\mathrm{Rest}$'' denotes a polynomial that does not contain the monomial $u_{s}^{j2^{s-1}}u_{s-1}^{j2^{s-2}}\cdots u_{2}^{2j} u_{1}^{j}$.
Hence, 
\begin{equation}
\label{eq:last}
q^j\notin \big\langle u_1^{j+1},u_2^{2j+1},u_3^{4j+1}\ldots, u_s^{j2^{s-1}+1}\big\rangle,
\end{equation}
in contradiction to relation~\eqref{eq:upper_bound}.
 
\subsection{}\label{subsec:compactness_argument}
We just proved that there cannot be any continuous $(\Z/2)^s$-equivariant map $Y \longrightarrow S(U_{s}^{\oplus j})$. 
Therefore, every continuous $(\Z/2)^s$-equivariant map $Y\longrightarrow U_{s}^{\oplus j}$ has a zero.

More is true: Since $Y$ is compact, every $(\Z/2)^s$-equivariant map $Y \longrightarrow U_{s}^{\oplus j}$, not necessarily continuous, which is a limit of a sequence of continuous $(\Z/2)^s$-equivariant maps $Y \longrightarrow U_{s}^{\oplus j}$, will also have a zero.
Indeed, let $\Psi:=\lim_{i\rightarrow\infty }\Psi^{i}$ where $\Phi^{i}\colon Y\longrightarrow U_{s}^{\oplus j}$ are continuous $(\Z/2)^s$-equivariant maps. 
Since the maps $\Psi^{i}$ are continuous, for every $i$ there exists a $y_i\in Y$ such that $\Psi^{i}(y_i)=0$.
The compactness of $Y$ yields the existence of a converging subsequence $\lim_{i\to\infty} y_{k(i)}= y\in Y$, where $k\colon\N\longrightarrow\N$ is a strictly increasing function.
Thus $\Psi(y)=\lim_{i\rightarrow\infty}\Psi^{k(i)}(y_{k(i)})=0$.

We have proved in \eqref{subsec : 2.5} that the maps $\phi_{\alpha,\gamma}$, used for the definition of the map $\Phi$, are limits of sequences of continuous maps.
Thus the $(\Z/2)^s$-equivariant map $\Phi\colon Y \longrightarrow U_{s}^{\oplus j}$, even not continuous by construction, still has a zero.
The proof of Theorem~\ref{th : main} is now complete.

\subsubsection*{Remark}
The choice of the degree bounds $\delta_\ell$ for the polynomials $p_\ell$ used for partitioning,
and consequently of the vector spaces $V_\ell$ etc., which in the special case $j=1$ already appeared in
Guth's work~\cite{guth2015}, can now be seen as very natural if one tries to 
show that at least one monomial in the power of the Dickson polynomial~$q$
does not belong to the index of the configuration space $Y$ of $s$-tuples of polynomials $(p_1,\dots,p_s)$,
and thus to obtain a contradiction in~\eqref{eq:last}.

\subsubsection*{Acknowledgements}
We are grateful to Josh Zahl and to the referee of JFPTA for very valuble remarks.


\end{document}